\documentclass[letterpaper, 10 pt, conference]{ieeeconf}  

\IEEEoverridecommandlockouts                              

\usepackage[normalem]{ulem}

\usepackage{amsmath,color}
\usepackage{amssymb}  
\usepackage{amsfonts}
\usepackage{graphicx}
\usepackage{dsfont}
\usepackage{psfrag}
\usepackage{epstopdf}
\usepackage{mathtools}
\usepackage{cite}
\usepackage{bbm}
\usepackage{algorithm}
\usepackage{algorithmic}
\usepackage{savesym}

\usepackage{amsthm}
 
\savesymbol{AND}

\theoremstyle{definition}
\newtheorem{definition}{Definition}
\newtheorem{assumption}{Assumption}
\newtheorem{remark}{Remark}

\newtheorem*{proofbf}{Proof}
\theoremstyle{plain}
\newtheorem{theorem}{Theorem}

\newtheorem{proposition}{Proposition}
\newtheorem{lemma}{Lemma}

\DeclareMathOperator*{\argmax}{argmax}

\newcommand{\EE}{\mathbb{E}}

\newcommand*\dd{\mathop{}\!\mathrm{d}}

\newcommand{\eqdef}{\overset{\textnormal{def}}{=}}

\usepackage{soul}

\title{\LARGE \bf
Convergence of the Expectation-Maximization Algorithm\\ Through Discrete-Time Lyapunov Stability Theory
}

\author{Orlando Romero$^{\dag}$\qquad Sarthak Chatterjee$^{\ddag}$\qquad S\'{e}rgio Pequito$^{\dag}$
\thanks{$^{\dag}$Department of Industrial and Systems Engineering, Rensselaer Polytechnic Institute, Troy NY, 12180, USA.}%
\thanks{$^{\ddag}$Department of Electrical, Computer, and Systems Engineering, Rensselaer Polytechnic Institute, Troy NY, 12180, USA.}%
}

\begin{document}

\maketitle
\thispagestyle{empty}
\pagestyle{empty}

\begin{abstract}
In this paper, we propose a dynamical systems perspective of the Expectation-Maximization (EM) algorithm. More precisely, we can analyze the EM algorithm as a nonlinear state-space dynamical system. The EM algorithm is widely adopted for data clustering and density estimation in statistics, control systems, and machine learning. This algorithm belongs to a large class of iterative algorithms known as \emph{proximal point methods}. In particular, we re-interpret limit points of the EM algorithm and other local maximizers of the likelihood function it seeks to optimize as equilibria in its dynamical system representation. Furthermore, we propose to assess its convergence as asymptotic stability in the sense of Lyapunov. As a consequence, we proceed by leveraging recent results regarding discrete-time Lyapunov stability theory in order to establish asymptotic stability (and thus, convergence) in the dynamical system representation of the EM algorithm. 

\end{abstract}

\section{Introduction}

With the ever-expanding size and complexity of data-sets used in the field of statistics, control systems, and machine learning, there has been a growing interest in developing algorithms that efficiently find the solution to the optimization problems that arise in these settings. For example, a fundamental problem in exploratory data mining is the problem of \emph{cluster analysis}, where the central task is to group objects into subgroups (\emph{i.e.}, \emph{clusters}) such that the objects in a particular cluster share several characteristics (or, features) with those that are, in some sense, sufficiently different from objects in different clusters \cite{Bishop2006,Tan:2005:IDM:1095618}.
\par The Expectation-Maximization (EM) algorithm \cite{Dempster1977maximum} is one of the most popular methods used in distribution-based clustering analysis and density estimation\cite{1212679,Mitchell:1997:ML:541177}. Given a dataset, we can assume that the data is distributed according to a finite mixture of Gaussian distributions whose parameters are randomly initialized and iteratively improved using the EM algorithm that seeks to maximize the likelihood that the data is justified by the distributions. This leads to finding the finite Gaussian mixture that hopefully best fits the dataset in question.

\par A current trend in optimization, machine learning, and control systems, is that of leveraging on a dynamical systems interpretation of iterative optimization algorithms \cite{LessardIQC,Elia2011,Mallada2018}.  The key idea is to view the the estimates themselves in the iterations of the algorithm as a state vector at different discrete instances of time (in particular, the initial approximation is viewed as the initial state), while the mechanism itself used to construct each subsequent estimate is modeled as a state-space dynamical system. Then, local optimizers and convergence in the optimization algorithm roughly translate to equilibria and asymptotic stability (in the sense of Lyapunov) in its dynamical system interpretation.

\par The convergence of the EM algorithm has been studied from the point-of-view of general point-to-set notions of convergence of optimization algorithms such as Zangwill's convergence theorem \cite{zangwill1969nonlinear}. Works such as \cite{wu1983convergence} and \cite{redner1984mixture} provide proofs of the convergence of the sequence of  estimates generated by the EM algorithm. 

The main contribution of this paper is to present a dynamical systems perspective of the convergence of the EM algorithm. The convergence of the EM algorithm is well known. However, our nonlinear stability analysis approach is intended to help open the field to new iteration schemes by possible addition of an artificial external input in the dynamical system representation of the EM algorithm. Then, leveraging tools from feedback systems theory, we could design a control law that translates to an accelerated convergence of the algorithm for specific subclassess of distributions.

The rest of the paper is organized as follows. In Section~\ref{sec:EMalg}, we briefly review the problem of maximum likelihood estimation and the EM algorithm. In Section~\ref{sec:dynamicalsysteminterpreation}, we propose a dynamical systems perspective of the EM algorithm and propose a particular generalized EM (GEM) algorithm. In Section~\ref{sec:localconvergence} we establish our main convergence results by leveraging discrete-time Lyapunov stability theory. Finally, Section~\ref{sec:conclusion} concludes the paper.

\subsection*{Notation}
The set of non-negative integers is represented by $\mathbb{Z}_+ = \{0,1,2,\ldots\}$, the set of real numbers is represented by $\mathbb{R}$, and $\mathbb{R}^n$ denotes the $n$-dimensional real vectors. The Euclidean norm is denoted by $\|\cdot\|$. We denote the open $\delta$-ball around a point $x\in\mathbb{R}^n$ as $B_\delta(x) = \{y\in\mathbb{R}^n: \|y-x\|<\delta\}$, and the closed $\delta$-ball as $\bar{B}_\delta(x) = \{y\in\mathbb{R}^n: \|y-x\|\leq\delta\}$. The gradient and Hessian matrix of a scalar function $f$ are denoted, respectively, by $\nabla f$ and $\nabla^2 f$. The notation $A\prec 0$ denotes that the real-valued square matrix $A$ is negative definite. We do not distinguish random vectors and their corresponding realizations through notation, but instead let it be implicit through context. For a given random vector $x\in\mathbb{R}^n$, we denote its probability distribution as $p(x)$. For simplicity, we will assume that every random vector is continuous, and hence, every distribution a probability density function. We denote the expected value of a function $f(x)$ of $x$ with respect to the distribution $p(x)$ by $\mathbb{E}_{p(x)}[f(x)] = \int p(x)f(x)\dd x$, or  $\mathbb{E}[f(x)]$ when the distribution $p(x)$ is clear from context.

\section{Expectation-Maximization Algorithm}
\label{sec:EMalg}

In this section, we recall the \emph{Expectation-Maximization} (EM) algorithm. Let $\theta\in\Theta\subseteq\mathbb{R}^p$ be some vector of unknown (but deterministic) parameters characterizing a distribution of interest, which we seek to infer from a collected dataset $y\in\mathbb{R}^m$ (from now on assumed fixed). 
To estimate $\theta$ from the dataset $y$, we first need a statistical model, \emph{i.e.}, an indexed class of probability distributions $\{p_\theta(y):\theta\in\Theta\}$. The function $\mathcal{L}:\Theta\to\mathbb{R}$ given by
\begin{equation}
    \mathcal{L}(\theta) = p_\theta(y)
    \label{eq:likelihoodfunction}
\end{equation}
denotes the \emph{likelihood function}.  The objective is to compute the \emph{maximum likelihood estimate} (MLE):
\begin{equation}
    \hat{\theta}_\textnormal{MLE} \eqdef \argmax_{\theta\in\Theta} \mathcal{L}(\theta),
    \label{ML}
\end{equation}
where the maximizer of $\mathcal{L}(\theta)$ is not necessarily unique, and hence, neither is the MLE. For that reason, it is actually more accurate to use
\begin{equation}
    \hat{\theta}_\textnormal{MLE} \in \argmax_{\theta\in\Theta} \mathcal{L}(\theta)
\end{equation}
as the definition of the MLE. From this point on, we will treat
\begin{equation}
    \argmax_{\theta\in\Theta}\mathcal{L}(\theta) = \left\{\theta\in\Theta: \mathcal{L}(\theta) = \max_{\theta'\in\Theta} \mathcal{L}(\theta')\right\}
\end{equation}
as a set, unless it consists of a single point $\theta^\star$, in which case we may use $\theta^\star = \argmax_{\theta\in\Theta} \mathcal{L}(\theta)$.

Next, we will introduce some assumptions that will ensure well-definedness throughout this paper.

\begin{assumption}
$\mathcal{L}(\theta) > 0$ for every $\theta\in\Theta$. \hfill $\circ$
\label{ass:supporty}
\end{assumption}
This assumption is simply a mild technical condition intended to avoid pathological behaviors, and is satisfied by most mixtures of distributions used in practice (\emph{e.g.}, Gaussian, Poisson, Beta). Furthermore, we surely have $\mathcal{L}(\theta)>0$ for at least \emph{some} $\theta\in\Theta$ (since, otherwise, the dataset $y$ is entirely useless regarding maximum likelihood estimation), and thus it suffices that we disregard from $\Theta$ any $\theta$ such that $\mathcal{L}(\theta) = 0$.

The underlying assumption for the EM algorithm is that there exists some latent (non-observable) random vector \mbox{$x\in\mathcal{X}\subseteq\mathbb{R}^n$} for which we possess a ``complete'' statistical model $\{p_\theta(x,y):x\in\mathcal{X},\theta\in\Theta\}$ (as opposed to the ``incomplete'' model $\{p_\theta(y):\theta\in\Theta\}$), and for which maximizing the expected value of the \emph{complete log-likelihood function} is easier than the \emph{incomplete likelihood function}. However, since $x$ is latent, the idea behind the EM algorithm is to iteratively maximize the expected complete log-likelihood. 

\begin{assumption}
$\mathcal{X} = \{x\in\mathbb{R}^n: p_\theta(x,y) > 0\}$ does not depend on $\theta\in\Theta$. \hfill $\circ$
\label{ass:supportx}
\end{assumption}

Together with Assumption~\ref{ass:supporty}, this assumption will further allow us to avoid certain pathological cases. Specifically, we can properly define the expected log-likelihood function (hereafter, also referred to as the \mbox{$Q$-function})~$Q:\Theta\times\Theta\to\mathbb{R}$, defined as
\begin{subequations}
\begin{align}
    Q(\theta,\theta') &\eqdef \EE_{p_{\theta'}(x|y)}[\log p_\theta(x,y)]\label{eq:Qfunction}\\
    &= \int_\mathcal{X} p_{\theta'}(x|y)\log p_\theta(x,y)\dd x,\label{eq:Qintegral}
\end{align}
\end{subequations}
where $p_{\theta'}(x|y) = p_{\theta'}(x,y)/p_{\theta'}(y)$. With all these ingredients and assumptions, we summarize the EM algorithm in Algorithm~\ref{alg:EM}. Notice that, the term $Q(\cdot,\theta_k)$ in Algorithm~\ref{alg:EM} denotes the expected complete log-likelihood function for any given iteration $k$. 
\begin{algorithm}
\caption{Expectation-Maximization (EM)}
\textbf{Input:} Observed data $y\in\mathbb{R}^m$, complete statistical model $\{p_\theta(x,y):x\in\mathcal{X},\theta\in\Theta\}$, and initial approximation $\theta_0$ of $\displaystyle\hat{\theta}_\textnormal{MLE}\in\argmax_{\theta\in\Theta}\mathcal{L}(\theta)$.\\
\textbf{Output:} $\theta_\infty$ such that hopefully $\displaystyle p_{\theta_\infty}(y) \approx \max_{\theta\in\Theta} \mathcal{L}(\theta)$.
\vspace*{-0.3cm}
\begin{algorithmic}[1]
\FOR{$k=0,1,2,\ldots$} 
    \STATE{\textbf{E-step:} \,\,Compute $Q(\theta,\theta_k)$ }
    \STATE{\textbf{M-step:} Determine $\displaystyle \theta_{k+1} \in \argmax_{\theta\in\Theta} Q(\theta,\theta_k)$}
\vspace{-0.15cm}
\ENDFOR
\RETURN{$\displaystyle\theta_\infty = \lim_{k\to\infty}\theta_k$, if it exists.}
\end{algorithmic}
\label{alg:EM}
\end{algorithm}

\begin{remark}
In practice, the iterations of Algorithm~\ref{alg:EM} are computed until some stopping criterion is achieved, such that it approximates $\theta_\infty$.\hfill $\circ$
\end{remark}

\section{Dynamical System Interpretation of the EM Algorithm and Convergence}
\label{sec:dynamicalsysteminterpreation}
Formally, the convergence of the EM algorithm is concerned with the existence and characteristics of the limit of the sequence $\{\theta_k\}_{k\in\mathbb{Z}_+}$ as $k\to\infty$. In particular, \emph{local convergence} refers to the property of the sequence $\{\theta_k\}_{k\in\mathbb{Z}_+}$ converging to the \emph{same} point $\theta^\star$ for every initial approximation $\theta_0$ that is sufficiently close to $\theta^\star$. On the other hand, \emph{global convergence} refers to convergence to the same point for \emph{any} initial approximation. Ideally, $\theta^\star$ is a global maximizer (or at least a local one) of the likelihood function. In practice, Algorithm~\ref{alg:EM} may converge to other stationary points of the likelihood function \cite{Dempster1977maximum,wu1983convergence}.

We will now see how the EM algorithm (such as many other iterative optimization algorithms) can be interpreted as a dynamical system in state-space, for which \emph{convergence} translates to \emph{(asymptotic) stability} in the sense of Lyapunov. 

To start, recall that any discrete-time time-invariant nonlinear dynamical system in state-space can be described by its dynamics, which are of the form
\begin{equation}
    \begin{cases}\theta[k+1] = F(\theta[k]),\qquad k\in\mathbb{Z}_+,\\ \theta[0] = \theta_0,\end{cases}
    \tag{$\mathcal{S}$}
    \label{eq:dynamics}
\end{equation}
and where $\theta[k]$ denotes the \emph{state} of the system and \mbox{$F:\Theta\to\Theta$} is some known function. In particular, any $F$ that satisfies
\begin{equation}
    F(\theta') \in \argmax_{\theta\in\Theta} Q(\theta,\theta')
\end{equation}
for every $\theta'\in\Theta$ represents a particular realization of the different iterations of Algorithm~\ref{alg:EM}. For the sake of simplicity, let us make the following assumption.
\begin{assumption}
$Q(\cdot,\theta')$ has a unique global maximizer in $\Theta$ for each $\theta'\in\Theta$.\hfill $\circ$
\label{ass:uniquemaximizerQ}
\end{assumption}

\begin{remark}
Assumption~\ref{ass:uniquemaximizerQ} does not imply that the likelihood function has a unique global maximizer, and subsequently, the MLE may still be non-unique. Furthermore, under Assumption~\ref{ass:uniquemaximizerQ}, the sequence $\{\theta_k\}_{k\in\mathbb{Z}_+}$ generated by Algorithm~\ref{alg:EM} is unique for each $\theta_0\in\Theta$, the function $F^\textnormal{EM}:\Theta\to\Theta$ given by
\begin{equation}
    F^{\textnormal{EM}}(\theta') = \argmax_{\theta\in\Theta} Q(\theta,\theta')
    \tag{\textsl{EM}}
    \label{eq:Ffunction}
\end{equation}
for $\theta'\in\Theta$ is uniquely defined, and \eqref{eq:dynamics} with $F=F^\textnormal{EM}$ captures the dynamical evolution emulated  by  Algorithm~\ref{alg:EM}, \emph{i.e.} $\theta[k] = \theta_k$ for every $k\in\mathbb{Z}_+$.\hfill $\circ$
\end{remark}

Recall that, for a dynamical system of the form~\eqref{eq:dynamics}, we say that $\theta^\star$ is an \emph{equilibrium} of the system if $\theta[0]=\theta^\star$ implies that $\theta[k] = \theta^\star$ for every $k\in\mathbb{Z}_+$. In other words, if $\theta^\star$ is a fixed point of $F(\theta)$, \emph{i.e.}, $F(\theta^\star)=\theta^\star$. For self-consistency, we now formally define Lyapunov stability.
\begin{definition}[Lyapunov stability]
Let $\theta^\star$ be an equilibrium of the dynamical system~\eqref{eq:dynamics}. We say that $\theta^\star$ is \emph{stable} if the trajectory $\theta[k]$ is arbitrarily close to $\theta^\star$ provided that it starts sufficiently close to $\theta^\star$. In other words, if, for any $\varepsilon > 0$, there exists some $\delta>0$ such that $\theta_0\in B_\delta(\theta^\star)$ implies that $\theta[k]\in B_\varepsilon(\theta^\star)$ for every $k\in\mathbb{Z}_+$.

Further, we say that $\theta^\star$ is \emph{(locally) asymptotically stable} if, apart from being stable, the trajectory $\theta[k]$ converges to $\theta^\star$ provided that it starts sufficiently close to $\theta^\star$. In other words, if $\theta^\star$ is stable and there exists some $\delta>0$ such that $\theta_0\in B_\delta(\theta^\star)$, implies that $\theta[k]\to\theta^\star$ as $k\to \infty$.\hfill $\circ$
\end{definition}

From the previous definition, it readily follows that asymptotic stability of the trajectory $\theta[k]$ generated by~\eqref{eq:dynamics} with $F=F^\textnormal{EM}$ is nearly equivalent to local convergence of the EM algorithm, since $\theta_k = \theta[k]\to\theta^\star$ for any $\theta_0$ in some sufficiently small open ball centered around $\theta^\star$. Therefore, to establish local convergence of the EM algorithm from the point of view of the asymptotic stability of the corresponding dynamical system, we first need to establish that the points of interest (\emph{i.e.}, the local maxima of the likelihood function) are equilibria of the system (\emph{i.e.}, fixed points of $F^\textnormal{EM}$).

Let $\theta^\star\in\Theta$ be a local maximizer of $\mathcal{L}(\theta)$. More precisely, let us start by considering $\theta^\star = \hat{\theta}_\textnormal{MLE} \in \argmax_{\theta\in\Theta}\mathcal{L}(\theta)$. Notice that the $Q$-function can be re-written as follows:
\begin{subequations}
\begin{align}
    Q(\theta,\theta') &= \EE_{p_{\theta'}(x|y)}[\log p_\theta(x,y)]  \\
    &= \log p_\theta(y) + \EE_{p_{\theta'}(x|y)}[\log p_\theta(x|y)] \\
    &= \log\mathcal{L}(\theta) - \mathcal{D}_\textnormal{KL}(\theta'\|\theta) - \mathcal{H}(\theta')\label{QShannon}
\end{align}
\end{subequations}
where $\mathcal{D}_\textnormal{KL}(\theta'\|\theta)$ denotes the Kullback-Leibler divergence from $p_\theta(\cdot|y)$ to $p_{\theta'}(\cdot|y)$, and $\mathcal{H}(\theta')$ denotes the differential Shannon entropy of $p_{\theta'}(\cdot|y)$ \cite{Cover:2006:EIT:1146355}.
Since the entropy term in~\eqref{QShannon} does not depend on  $\theta$, then
\begin{equation}
    F^\textnormal{EM}(\theta') = \argmax_{\theta\in\Theta} \,\{\log\mathcal{L}(\theta) - \mathcal{D}_\textnormal{KL}(\theta'\|\theta)\}
\label{eq:FShannon}
\end{equation}
for $\theta'\in\Theta$.

\begin{remark}
In general, the Kullback-Leibler divergence from an arbitrary distribution $q(x)$ to another $p(x)$, denoted as $\mathcal{D}_\textnormal{KL}(p\|q)$, satisfies the following two properties: \emph{(i)} $\mathcal{D}_\textnormal{KL}(p\|q) \geq 0$ for every $p$ and $q$, a result known as \emph{Gibbs' inequality}; and \emph{(ii)} $\mathcal{D}_\textnormal{KL}(p\|q) = 0$ if and only if $p=q$ almost everywhere. \hfill $\circ$
\end{remark}

Therefore, we can now state the following.

\begin{proposition}
The MLE is an equilibrium of the dynamical system \eqref{eq:dynamics} with $F=F^\textnormal{EM}$.\hfill $\diamond$
\label{thm:MLEequilibrium}
\end{proposition}

\begin{proofbf}
Inspecting \eqref{eq:FShannon} at $\theta' = \hat{\theta}_\textnormal{MLE} = \argmax_{\theta\in\Theta}\mathcal{L}(\theta)$, it readily follows that $\log\mathcal{L}(\theta)$ and $-\mathcal{D}_\textnormal{KL}(\hat{\theta}_\textnormal{MLE}\|\theta)$ are both maximized at $\theta = \hat{\theta}_\textnormal{MLE}$, which in turn implies that $F^\textnormal{EM}(\hat{\theta}_\textnormal{MLE}) = \hat{\theta}_\textnormal{MLE}$.\hfill$\blacksquare$
\end{proofbf}

However, it should be clear that the previous argument does not hold for non-global maximizers of the likelihood function. One approach to get around this issue is to consider a specific variant of the \emph{generalized} EM algorithm (GEM)\footnote{A GEM algorithm is any variant of the EM algorithm where the M-step is replaced by a search of some $\theta_{k+1}\in\Theta$ such that $Q(\theta_{k+1},\theta_k) > Q(\theta_k,\theta_k)$, if one exists (otherwise $\theta_{k+1} = \theta_k$), not necessarily in $\argmax_{\theta\in\Theta} Q(\theta,\theta_k)$.}. More precisely, we propose a GEM algorithm that searches for a global maximizer of $Q(\theta,\theta_k)$ in a restricted parameter space: $B_\delta(\theta_k)$ (\emph{i.e.}, a $\delta$-ball around $\theta_k$). We call this algorithm the $\delta$-EM algorithm, which is summarized in Algorithm~\ref{alg:deltaEM}.

\begin{algorithm}
\caption{$\delta$-Expectation-Maximization ($\delta$-EM)}
\textbf{Input:} Restricted parameter space radius $\delta>0$, observed data $y\in\mathbb{R}^m$, complete statistical model $\{p_\theta(x,y):x\in\mathcal{X},\theta\in\Theta\}$, and initial approximation $\theta_0$ of $\displaystyle\hat{\theta}_\textnormal{MLE}\in\argmax_{\theta\in\Theta}\mathcal{L}(\theta)$.\\
\textbf{Output:} $\theta_\infty$ such that hopefully $\displaystyle p_{\theta_\infty}(y) \approx \max_{\theta\in\Theta} p_\theta(y)$.
\vspace*{-0.3cm}
\begin{algorithmic}[1]
\FOR{$k=1,2,\ldots$} 
    \STATE{\textbf{E-step:} \,\,Compute $Q(\theta,\theta_k)$ }
    \STATE{\textbf{M-step:} Determine $\displaystyle \theta_{k+1} \in \argmax_{\theta\in B_\delta(\theta_k)\cap\Theta} Q(\theta,\theta_k)$}
\vspace{-0.15cm}
\ENDFOR
\RETURN{$\displaystyle\theta_\infty = \lim_{k\to\infty}\theta_k$, if it exists.}
\end{algorithmic}
\label{alg:deltaEM}
\end{algorithm}

\noindent Next, we make the following simplifying assumption.
\begin{assumption}
$Q(\cdot,\theta')$ has a unique global maximizer in $B_\delta(\theta')\cap\Theta$, for every $\theta'\in\Theta$ and $\delta>0$.\hfill $\circ$
\label{ass:uniquemaximizerQdelta}
\end{assumption}

\noindent Naturally, we can interpret Algorithm~\ref{alg:deltaEM} as the dynamical system \eqref{eq:dynamics} with $F=F^{\delta-\textnormal{EM}}:\Theta\to\Theta$ given by
\begin{align}
    F^{\delta-\textnormal{EM}}(\theta') &\eqdef \argmax_{\theta\in B_\delta(\theta_k)\cap\Theta} Q(\theta,\theta')\tag{$\delta$-\textsl{EM}}\label{eq:FdeltaEM}\\
    &= \argmax_{\theta\in B_\delta(\theta')\cap\Theta} \,\{\log\mathcal{L}(\theta) - \mathcal{D}_\textnormal{KL}(\theta'\|\theta)\},\label{eq:FdeltaShannon}
\end{align}
where \eqref{eq:FdeltaShannon} was derived following a similar argument that led to \eqref{eq:FShannon}. Therefore, similar to Proposition~1, we have the following result. 

\begin{proposition}
Any local maximizer of the likelihood function is an equilibrium of the dynamical system \eqref{eq:dynamics} with $F=F^{\delta-\textnormal{EM}}$, provided that $\delta>0$ is small enough.\hfill $\diamond$
\label{thm:deltaEMequilibrium}
\end{proposition}

\begin{proofbf}
Let $\theta^\star\in\Theta$ be a local maximizer of $\mathcal{L}(\theta)$. Inspecting \eqref{eq:FdeltaShannon} at $\theta'=\theta^\star$, it readily follows that $-\mathcal{D}_\textnormal{KL}(\theta^\star\|\theta)$ is maximized at $\theta=\theta^\star$. Furthermore, if $\delta>0$ is small enough, then $\log\mathcal{L}(\theta)$ is also maximized (in $B_\delta(\theta^\star)\cap\Theta$) at $\theta=\theta^\star$, which implies that $F^{\delta-\textnormal{EM}}(\theta^\star) = \theta^\star$.\hfill $\blacksquare$
\end{proofbf}

\section{Local Convergence Through Discrete-Time Lyapunov Stability Theory}
\label{sec:localconvergence}

In this section, we will discuss how we can establish the local convergence of the EM algorithm by exploiting classical results from Lyapunov stability theory in the dynamical system interpretation of the EM algorithm.

To start, we state the discrete-time version of the Lyapunov theorem (Theorem 1.2 in \cite{discretelyapunov}). First, recall that a function $\mathcal{V}:\Theta\to\mathbb{R}$ is said to be \emph{positive semidefinite}, if $\mathcal{V}(\theta)\geq 0$ for every $\theta\in\Theta$. Furthermore, we say that $\mathcal{V}$ is \emph{positive definite} with respect to $\theta^\star\in\Theta$, if $\mathcal{V}(\theta^\star) = 0$ and $\mathcal{V}(\theta)>0$ for $\theta\in\Theta\setminus\{\theta^\star\}$. 

\begin{theorem}[Lyapunov Stability]
Let $\theta^\star\in\Theta$ be an equilibrium of the dynamical system \eqref{eq:dynamics} in the interior of $\Theta$ and let $\delta>0$ be such that $B_\delta(\theta^\star)\subseteq\Theta$. Suppose that $F$ is continuous and there exists some continuous function $\mathcal{V}:B_\delta(\theta^\star)\to\mathbb{R}$ (called a \emph{Lyapunov function}) such that $\mathcal{V}$ and $-\Delta\mathcal{V}$ are, respectively, positive definite with respect to $\theta^\star$ and positive semidefinite, where $\Delta\mathcal{V}(\theta) \eqdef \mathcal{V}(F(\theta)) - \mathcal{V}(\theta)$. Then, $\theta^\star$ is stable. If $-\Delta\mathcal{V}$ is also positive definite with respect to $\theta^\star$, then $\theta^\star$ is asymptotically stable.\hfill $\diamond$
\label{thm:lyapunov}
\end{theorem}

\begin{remark}
The proof of Theorem~\ref{thm:lyapunov} can be found in \cite{discretelyapunov}. The statement of the theorem therein assumes local Lipschitz continuity of $F$. This is likely a residual from the classical assumption of local Lipschitz continuity of $F$ in continuous systems with dynamics of the form $\dot{\theta}(t) = F(\theta(t))$, which is required by the Picard-Lindel\"{o}f theorem to ensure unique existence of a solution to the differential equation $\dot{\theta}(t)=F(\theta(t))$ for each initial state $\theta(0)=\theta_0$. For discrete-time systems, on the other hand, the unique existence of the trajectory is immediate. However, a careful analysis of the argument used in the proof found in \cite{discretelyapunov} reveals that the continuity of $F$ is nevertheless implicitly needed to ensure the continuity of $\Delta\mathcal{V}$, since the extreme value theorem is invoked for $\Delta\mathcal{V}$.\hfill $\circ$
\end{remark}

In order to leverage the previous theorem to establish local convergence of the EM algorithm to local maxima of the likelihood function, we need to propose a \emph{candidate} Lyapunov function. However, before doing so, we need to ensure that $F=F^\textnormal{EM}$ is continuous, which is attained by imposing some regularity on the likelihood function. 

\begin{assumption}
$\mathcal{L}(\theta)$ is twice continuously differentiable. \hfill $\circ$
\label{ass:C2}
\end{assumption}

Subsequently, we obtain the following result.

\begin{lemma}
$F^\textnormal{EM}$ and $F^{\delta-\textnormal{EM}}$ are both continuous. \hfill $\diamond$
\label{lemma1}
\end{lemma}

\begin{proofbf}
Under Assumption~\ref{ass:C2}, it readily follows from \eqref{eq:Qintegral} that $Q(\cdot,\cdot)$ is continuous in both of its arguments. Let $\{\theta_k'\}_{k\in\mathbb{Z}_+}\subseteq\Theta$ be a sequence converging to $\theta'\in\Theta$. Note that, for each $k\in\mathbb{Z}_+$, we have $Q(F^\textnormal{EM}(\theta_k'),\theta_k') \geq Q(\theta,\theta_k')$ for every $\theta\in\Theta$. Taking the limit when $k\to\infty$, and leveraging the continuity of $Q$, we have $Q\left(\lim_{k\to\infty}F^\textnormal{EM}(\theta_k'),\theta'\right) \geq Q(\theta,\theta')$ for every $\theta\in\Theta$. Consequently,
\begin{equation}
    Q\left(\lim_{k\to\infty}F^\textnormal{EM}(\theta_k'),\theta'\right) = \max_{\theta\in\Theta}Q(\theta,\theta'),
\end{equation}
and therefore,
\begin{equation}
    \lim_{k\to\infty}F^\textnormal{EM}(\theta_k') = \argmax_{\theta\in\Theta} Q(\theta,\theta') = F^\textnormal{EM}(\theta').
\end{equation}
This same argument can be readily adapted for $F^{\delta-\textnormal{EM}}$. \hfill $\blacksquare$
\end{proofbf}


Let $\theta^\star\in\Theta$ be a local maximizer of $\mathcal{L}(\theta)$. Once again, let us start by considering $\theta^\star = \hat{\theta}_\textnormal{MLE}\in\argmax_{\theta\in\Theta}\mathcal{L}(\theta)$. From Proposition~\ref{thm:MLEequilibrium}, we know that $\theta^\star$ is an equilibrium of~\eqref{eq:dynamics} for $F=F^\textnormal{EM}$. Next, a naive guess of a candidate Lyapunov function would be to consider $\mathcal{V}(\theta) = \mathcal{L}(\theta)$, since this would satisfy $\mathcal{V}(\theta) \geq 0$ for every $\theta\in\Theta$, but  $\mathcal{V}(\theta^\star) > 0$; hence, it is not a Lyapunov function. Notwithstanding, if  we subtract $\mathcal{L}(\theta^\star)$ from the previous candidate, \emph{i.e.}, $\mathcal{V}(\theta) = \mathcal{L}(\theta) - \mathcal{L}(\theta^\star)$,
then $\mathcal{V}(\theta^\star) = 0$ and $\mathcal{V}(\theta) < 0$ for $\theta\in\Theta\setminus\argmax_{\theta'\in\Theta} \mathcal{L}(\theta')$. As a consequence, it should be clear that
\begin{equation}
    \mathcal{V}(\theta) = \mathcal{L}(\theta^\star) - \mathcal{L}(\theta)
    \label{eq:lyapfunc}
\end{equation}
appears to be the ideal candidate, since $\mathcal{V}(\theta)\geq 0$ for every $\theta\in\Theta$ and $\mathcal{V}(\theta^\star) = 0$. Yet, $\mathcal{V}$ may be only positive semidefinite instead of positive definite (with respect to $\theta^\star$), since $\mathcal{V}(\theta) = 0$ if and only if $\theta\in \argmax_{\theta'\in\Theta} \mathcal{L}(\theta')$. To circumvent this issue, we will need to assume that $\theta^\star$ is an isolated maximizer\footnote{We say that $\theta^\star\in\Theta$ is an \emph{isolated} maximizer (stationary point) of $\mathcal{L}(\theta)$ if it is the only local maximizer (stationary point) of $\mathcal{L}(\theta)$ in $B_r(\theta^\star)\cap\Theta$ for some small enough $r>0$.} of $\mathcal{L}(\theta)$. 

\begin{lemma}
Suppose that $\theta^\star\in\Theta$ is an isolated maximizer of $\mathcal{L}(\theta)$. Then, $\mathcal{V}:B_r(\theta^\star)\to\mathbb{R}$ given by \eqref{eq:lyapfunc} with $F=F^\textnormal{EM}$ or $F=F^{\delta-\textnormal{EM}}$ is positive definite with respect to $\theta^\star$ and $-\Delta\mathcal{V}$ is positive semidefinite, provided that $r>0$ is small enough. \hfill $\diamond$
\label{lemma:Visgood}
\end{lemma}

\begin{proofbf}
Let us focus on the case $F=F^\textnormal{EM}$. The positive definiteness follows from~\eqref{eq:lyapfunc} and the definition of isolated maximizer. 
On the other hand, from \eqref{eq:FShannon}, it follows that
\begin{equation}
    \begin{split}
        \log\mathcal{L}(F^\textnormal{EM}(\theta)) - &\mathcal{D}_\textnormal{KL}(\theta\|F^\textnormal{EM}(\theta))\nonumber\\
    &\geq \log\mathcal{L}(\theta) - \underbrace{\mathcal{D}_\textnormal{KL}(\theta\|\theta)}_{=0},
    \end{split}
    \label{eq:abc}
\end{equation}
and therefore,
\begin{equation}
    \log\mathcal{L}(F^\textnormal{EM}(\theta)) \geq \log\mathcal{L}(\theta) + \underbrace{\mathcal{D}_\textnormal{KL}(\theta\,\|\,F^\textnormal{EM}(\theta))}_{\geq 0}\\
\label{eq:deltaVnonnegative}
\end{equation}
for every $\theta\in\Theta$. Thus, from the strict monotonicity of the logarithm function, it follows that $-\Delta\mathcal{V}(\theta) = \mathcal{L}(F^\textnormal{EM}(\theta))-\mathcal{L}(\theta)$ is indeed positive semidefinite. The case $F=F^{\delta-\textnormal{EM}}$ follows by essentially the same argument.  \hfill $\blacksquare$
\end{proofbf}

\begin{remark}
We are now in conditions to establish the stability of isolated MLEs in the dynamical system that represents the EM algorithm. Further, through a similar argument, the stability of arbitrary isolated local maximizers of the likelihood for the dynamical system that represents the $\delta$-EM algorithm with small enough $\delta>0$. However, non-asymptotic stability does not seem to translate into any interesting aspect of the convergence of the EM or $\delta$-EM algorithms.\hfill $\circ$
\end{remark}


In order to establish the positive definiteness of $-\Delta\mathcal{V}$, we first need to characterize the equilibria of the dynamical systems that represent Algorithms~\ref{alg:EM} and \ref{alg:deltaEM}. 


\begin{lemma}
Every fixed point $F^\textnormal{EM}$ and $F^{\delta-\textnormal{EM}}$ in the interior of $\Theta$ is a stationary point of the likelihood function. \hfill $\diamond$
\label{lemma2}
\end{lemma}

\begin{proofbf}
Let $\theta^\star\in\Theta$ be a fixed point of $F^\textnormal{EM}$ or $F^{\delta-\textnormal{EM}}$ in the interior of $\Theta$. From Assumption \ref{ass:C2}, and from \eqref{eq:FShannon} and \eqref{eq:FdeltaShannon}, it follows that
\begin{equation}
    \frac{\partial}{\partial\theta}\Big\{\log \mathcal{L}(\theta) - \mathcal{D}_\textnormal{KL}(\theta^\star\|\theta) \Big\}\Big\vert_{\theta=\theta^\star} = 0,
\end{equation}
where the divergence term actually vanishes, since $\theta=\theta^\star$ is a (global) minimizer of $\mathcal{D}_\textnormal{KL}(\theta^\star\|\theta)$. It thus follows that $\theta^\star$ is a stationary point of $\mathcal{L}(\theta)$. \hfill $\blacksquare$
\end{proofbf}

Next, we will need to assume additional regularity on the likelihood function. A very common assumption for a parameterized statistical model is for the parameterization to be injective, meaning that each distribution in the model is indexed by exactly one instance of the parameter space.

\begin{assumption}
$\theta\mapsto p_\theta(\cdot|y)$ is injective. \hfill $\circ$
\label{ass:injectiveparameterization}
\end{assumption}


Equipped with the last two assumptions (Assumption~\ref{ass:C2} and \ref{ass:injectiveparameterization}), we are ready to establish the positive definiteness of $-\Delta\mathcal{V}$, and subsequently, the local convergence of Algorithm~\ref{alg:EM} to ML estimates.

\begin{theorem}[Local Convergence of EM to MLE]
If $\nabla^2\mathcal{L}(\hat{\theta}_\textnormal{MLE})\prec 0$, then the sequence $\{\theta_k\}_{k\in\mathbb{Z}_+}$ generated by Algorithm~\ref{alg:EM} converges to $\hat{\theta}_\textnormal{MLE} = \argmax_{\theta\in\Theta} \mathcal{L}(\theta)$ for every initial approximation $\theta_0\in\Theta$ in a small enough open ball centered around $\hat{\theta}_\textnormal{MLE}$. \hfill $\diamond$
\label{thm:localconvergenceMLE}
\end{theorem}


\begin{proofbf}
First, recall from Proposition~\ref{thm:MLEequilibrium} that $\hat{\theta}_\textnormal{MLE}$ is an equilibrium of \eqref{eq:dynamics} with $F=F^\textnormal{EM}$ (\emph{i.e.}, a fixed point of $F^\textnormal{EM}$), and from Lemma~\ref{lemma1} that $F^\textnormal{EM}$ is continuous. Furthermore, $\hat{\theta}_\textnormal{MLE}$ is in the interior of $\Theta$ since $\nabla^2\mathcal{L}(\hat{\theta}_\textnormal{MLE})\prec 0$.

Let $r>0$ be small enough such that $\hat{\theta}_\textnormal{MLE}$ is the only stationary point of $\mathcal{L}(\theta)$ in $B_r(\hat{\theta}_\textnormal{MLE})\cap\Theta$. Such $r>0$ can be chosen since $\hat{\theta}_\textnormal{MLE}$ is itself a stationary point of $\mathcal{L}(\theta)$ and $\nabla^2\mathcal{L}(\hat{\theta}_\textnormal{MLE})\prec 0$ (which, together, they ensure isolated stationarity). In particular, $\hat{\theta}_\textnormal{MLE}$ is an isolated maximizer. Then, from Lemma~\ref{lemma:Visgood}, the function $\mathcal{V}:B_r(\hat{\theta}_\textnormal{MLE})\to\mathbb{R}$ given by \eqref{eq:lyapfunc} with $\theta^\star = \hat{\theta}_\textnormal{MLE}$ is positive definite with respect to $\hat{\theta}_\textnormal{MLE}$, and $-\Delta\mathcal{V}$ is positive semidefinite. 

Note that the inequality of the divergence term in \eqref{eq:deltaVnonnegative} is strict for $\theta\in B_r(\hat{\theta}_\textnormal{MLE})\setminus\{\hat{\theta}_\textnormal{MLE}\}$. This is because, from Assumption~\ref{ass:injectiveparameterization}, $\mathcal{D}_\textnormal{KL}(\theta\|F^\textnormal{EM}(\theta)) = 0$ if and only if $\theta$ is a fixed point of $F^\textnormal{EM}$. But such a point needs to be a stationary point (see Lemma~\ref{lemma2}), which would lead to the contradiction $\theta=\hat{\theta}_\textnormal{MLE}$, since $\hat{\theta}_\textnormal{MLE}$ is the only stationary point of $\mathcal{L}(\theta)$ in $B_r(\hat{\theta}_\textnormal{MLE})\setminus\{\hat{\theta}_\textnormal{MLE}\}$. Therefore, $\log\mathcal{L}(F^\textnormal{EM}(\theta)) > \log\mathcal{L}(\theta)$, \emph{i.e.}, $-\Delta\mathcal{V}(\theta) = \mathcal{L}(F^\textnormal{EM}(\hat{\theta}_\textnormal{MLE})) - \mathcal{L}(\theta) > 0$ for every $\theta\in B_r(\hat{\theta}_\textnormal{MLE})\setminus\{\hat{\theta}_\textnormal{MLE}\}$. Finally, since $\hat{\theta}_\textnormal{MLE}$ is a fixed point of $F^\textnormal{EM}$, then $\Delta\mathcal{V}(\hat{\theta}_\textnormal{MLE}) = 0$, which concludes that $-\Delta\mathcal{V}$ is positive definite. The conclusion follows by invoking Theorem~\ref{thm:lyapunov}, since $\hat{\theta}_\textnormal{MLE}$ was just proved to be asymptotically stable. \hfill $\blacksquare$
\end{proofbf}


\begin{theorem}[Local Convergence of $\delta$-EM to Local Maxima]
If $\theta^\star\in\Theta$ is such that $\nabla\mathcal{L}(\hat{\theta}^\star) = 0$ and $\nabla^2\mathcal{L}(\theta^\star)\prec 0$, and $\delta>0$ is small enough, then the sequence $\{\theta_k\}_{k\in\mathbb{Z}_+}$ generated by Algorithm~\ref{alg:deltaEM} converges to $\theta^\star$ for every initial approximation $\theta_0\in\Theta$ in a small enough open ball centered around $\theta^\star$. \hfill $\diamond$
\label{thm:convergencedeltaEM}
\end{theorem}

\begin{proofbf}
The proof follows similar steps to those in the proof of Theorem ~\ref{thm:localconvergenceMLE} with the following adaptations: first, replace $F^\textnormal{EM}$ by $F^{\delta-\textnormal{EM}}$, and secondly,  $\hat{\theta}_\textnormal{MLE}$ by $\theta^\star$. \hfill $\blacksquare$
\end{proofbf}

Notice that, the reason why Theorem~\ref{thm:convergencedeltaEM} cannot be readily adapted for the EM algorithm (as opposed to the $\delta$-EM algorithm) is that local maximizers may fail to be fixed points of $F^\textnormal{EM}$ and therefore, equilibria of \eqref{eq:dynamics} with $F=F^\textnormal{EM}$. To circumvent this limitation, we will focus on the \emph{limit points} of the EM algorithm. First, recall that $\theta^\star\in\Theta$ is a fixed point of Algorithm~\ref{alg:EM}, if there exists some $\theta_0\in\Theta$ such the $\theta_k\to \theta^\star$ as $k\to\infty$ for the sequence $\{\theta_k\}_{k\in\mathbb{Z}_+}$ generated by Algorithm~\ref{alg:EM}, which is captured by the following result.

\begin{lemma}
If $\theta^\star\in\Theta$ is a limit point of Algorithm~\ref{alg:EM}, then it is also a fixed point of $F^\textnormal{EM}$. \hfill $\diamond$
\label{lemma:limitpoint}
\end{lemma}

\begin{proofbf}
Let $\theta_0\in\Theta$ be such that $\theta_k\to\theta^\star$ as $k\to\infty$, where $\{\theta_k\}_{k\in\mathbb{Z}_+}$ was generated by Algorithm~\ref{alg:EM}. Then, by the continuity of $F^\textnormal{EM}$, it follows that $F^\textnormal{EM}(\theta^\star) = \lim_{k\to\infty} F^\textnormal{EM}(\theta_k) = \lim_{k\to\infty} \theta_{k+1} = \theta^\star$. \hfill $\blacksquare$
\end{proofbf}

\begin{remark}
Notice that, while not every local maximizer of the likelihood function is a limit point of the EM algorithm, the same is not true for the $\delta$-EM algorithm, provided that $\delta>0$ is sufficiently small and the local maximizer is sufficiently regular (\emph{i.e.}, an isolated stationary point of the likelihood function).
\label{remark:limitpoints}
\end{remark}

Upon Remark~\ref{remark:limitpoints}, and the convergence results established before, we can now establish the following claim.

\begin{theorem}[Local Convergence of EM to its Limit Points]
If $\theta^\star\in\Theta$ is a limit point of Algorithm~\ref{alg:EM} such that $\nabla^2\mathcal{L}(\theta^\star)\prec 0$, then the sequence $\{\theta_k\}_{k\in\mathbb{Z}_+}$ generated by Algorithm~\ref{alg:EM} converges to $\theta^\star$ for every initial approximation $\theta_0\in\Theta$ in a small enough open ball centered around $\theta^\star$. \hfill $\diamond$
\label{thm:localconvergenceEMlimitpoint}
\end{theorem}

\begin{proofbf}
The proof follows similar steps to those in the proof of Theorem ~\ref{thm:localconvergenceMLE}, where $\hat{\theta}_\textnormal{MLE}$ is replaced by $\theta^\star$, and followed by invoking Lemma~\ref{lemma:limitpoint} instead of Proposition~\ref{thm:MLEequilibrium}. \hfill $\blacksquare$
\end{proofbf}

We conclude this section by exploring how the notion of \emph{exponential stability} can be leveraged to bound the \emph{convergence rate} of the EM algorithm. First, recall that the (linear) rate of convergence for a sequence $\{\theta_k\}_{k\in\mathbb{Z}_+}$ is the number $0\leq \mu \leq 1$ given by
\begin{equation}
    \mu = \lim_{k\to\infty} \frac{\|\theta_{k+1} - \theta^\star\|}{\|\theta_k - \theta^\star\|},
\end{equation}
provided that the limit exists. Additionally, recall that an equilibrium $\theta^\star\in\Theta$ of \eqref{eq:dynamics} is said to be \emph{exponentially stable} if there exist constants $c,\gamma >0$ such that, for every $\theta_0\in\Theta$ in a sufficiently small open ball centered around $\theta^\star$, we have  $\|\theta[k] - \theta^\star\| \leq c\cdot e^{-\gamma k}\|\theta_0 - \theta^\star\|$ for every $k\in\mathbb{Z}_+$.
\begin{remark}
If $\theta^\star$ is an exponentially stable equilibrium of \eqref{eq:dynamics} with $F=F^\textnormal{EM}$, then $\{\theta_k\}_{k\in\mathbb{Z}}$ generated by Algorithm~\ref{alg:EM} converges to $\theta^\star$ with linear convergence rate 
\begin{subequations}
\begin{align}
\mu &= \lim_{k\to\infty} \frac{\|\theta_{k+1} - \theta^\star\|}{\|\theta_k - \theta^\star\|}\\
&\leq \lim_{k\to\infty} \frac{c\cdot e^{-\gamma\times 0}\|\theta_k - \theta^\star\|}{\|\theta_k - \theta^\star\|}\\
&= c,
\end{align}
\end{subequations}
for every initial approximation $\theta_0\in\Theta$ in a sufficiently small open ball centered around $\theta^\star$. \hfill $\circ$
\end{remark}

The following theorem (adapted from Theorem 5.7 in \cite{discretelyapunov}) allows us to ensure exponential stability, and subsequently to bound the linear convergence rate, provided our Lyapunov function satisfies some additional (mild) regularity conditions.

\begin{theorem}[Exponential Stability]
Let $\theta^\star\in\Theta$ be an equilibrium of \eqref{eq:dynamics} in the interior of $\Theta$, with $B_\delta(\theta^\star)\subseteq\Theta$ for some small enough $\delta>0$, and $F:\Theta\to\Theta$ continuous. Let $\mathcal{V}:B_\delta(\theta^\star)\to\mathbb{R}$ be a continuous and positive definite function (with respect to $\theta^\star$) such that
\begin{subequations}
\begin{align}
    \mathcal{V}(\theta) &\leq a\|\theta - \theta^\star\|^2,\label{eq:Vexp1}\\
    -\Delta\mathcal{V}(\theta) &\geq b\|\theta-\theta^\star\|^2,
    \label{eq:Vexp2}
\end{align}
\label{eq:Vexp}
\end{subequations}
for every $\theta\in\Theta$, for some constants $a,b > 0$. Then, $\theta^\star$ is exponentially stable. More precisely, we have $\|\theta[k]-\theta^\star\|\leq c\cdot e^{-\gamma k}$ for every $k\in\mathbb{Z}$ and $\theta_0\in B_r(\theta^\star)$ with $r>0$ small enough, where $c=d/a$ with
\begin{equation}
    d = \lim_{\delta\to 0} \max_{\theta\in \bar{B}_r(\theta^\star)\setminus B_\delta(\theta^\star)} \frac{\mathcal{V}(\theta)}{\|\theta-\theta^\star\|},\label{eq:d}
\end{equation}
and $\gamma = \log a-\log(a-b)$.\hfill $\diamond$
\label{thm:exponentialLyapunov}
\end{theorem}

Equipped with this result, we are now ready to establish the following sufficient conditions.

\begin{proposition}
Let $\theta^\star\in\Theta$ be a limit point of Algorithm~\ref{alg:EM} such that $\nabla^2\mathcal{L}(\theta^\star) \prec 0$. Suppose that there exist constants $a,b>0$ such that
\begin{subequations}
\begin{align}
\mathcal{L}(\theta) - \mathcal{L}(\theta^\star) &\geq - a\|\theta - \theta^\star\|^2\label{eq:LminusLa},\\
\mathcal{L}(F^\textnormal{EM}(\theta)) - \mathcal{L}(\theta) &\geq b\|\theta-\theta^\star\|^2,\label{eq:LminusLb}
\end{align}
\label{eq:LminusL}
\end{subequations}
or
\begin{subequations}
\begin{align}
\mathcal{L}(\theta)/\mathcal{L}(\theta^\star) &\geq \exp \{- a\|\theta - \theta^\star\|^2\},\\
\mathcal{D}_\textnormal{KL}(\theta\|F^\textnormal{EM}(\theta)) &\geq b\|\theta-\theta^\star\|^2,. 
\end{align}
\label{eq:LfracL}
\end{subequations}
for every $\theta\in B_\delta(\theta^\star)$ in some small enough $\delta>0$. Then, $\theta^\star$ is exponentially stable.\hfill $\diamond$
\label{prop:exp}
\end{proposition}
\begin{proofbf}
From Lemma~\ref{lemma:limitpoint}, it follows that $\theta^\star$ is an equilibrium of \eqref{eq:dynamics} with $F=F^\textnormal{EM}$. The result follows from invoking Theorem~\ref{thm:exponentialLyapunov}. First, to see that condition~\eqref{eq:Vexp1} is verified, we notice that this is equivalent to \eqref{eq:LminusLa} for $\mathcal{V}(\theta) = \mathcal{L}(\theta^\star) - \mathcal{L}(\theta)$ (which is continuous and positive definite with respect to $\theta^\star$ in $B_\delta(\theta^\star)$). On the other hand, condition~\eqref{eq:Vexp2} readily follows from~\eqref{eq:LminusLb}. Similarly, \eqref{eq:Vexp} follows from \eqref{eq:LfracL} for $\mathcal{V}(\theta) = \log\mathcal{L}(\theta^\star) - \log\mathcal{L}(\theta)$ (which is also continuous and positive definite with respect to $\theta^\star$ in $B_\delta(\theta^\star)$). \hfill $\blacksquare$
\end{proofbf}

\begin{remark}
Notice that \eqref{eq:LminusLa} actually readily holds for $\displaystyle a = -\frac{1}{2}\min_{\theta\in B_\delta(\theta^\star)} \lambda_\textnormal{min}[\nabla^2\mathcal{L}(\theta)]$. \hfill $\circ$
\end{remark}

Lastly, as consequence of Theorems~\ref{thm:localconvergenceEMlimitpoint} and \ref{thm:exponentialLyapunov}, we have the following result. 
\begin{theorem}[Explicit Bound for Convergence Rate of EM]
Under the conditions of Proposition~\ref{prop:exp}, the linear convergence rate of Algorithm~\ref{alg:EM} can be bounded as $\mu\leq d/a$, where \eqref{eq:d} defined through $\mathcal{V}(\theta) = \mathcal{L}(\theta^\star) - \mathcal{L}(\theta)$ for the case \eqref{eq:LminusL}, and $\mathcal{V}(\theta) = \log\mathcal{L}(\theta^\star) - \log\mathcal{L}(\theta)$ for the case \eqref{eq:LfracL}. \hfill $\diamond$
\end{theorem}

\section{Conclusion and Future Work}
\label{sec:conclusion}
We proposed a dynamical systems perspective of the Expectation-Maximization (EM) algorithm, by analyzing the evolution of its estimates as a nonlinear state-space dynamical system. In particular, we drew parallels between limit points and equilibria, convergence and asymptotic stability, and we leveraged results on discrete-time Lyapunov stability theory to establish local convergence results for the EM algorithm. In particular, we derived conditions that allow us to construct explicit bounds on the linear convergence rate of the EM algorithm by establishing exponential stability in the dynamical system that represents it. Future work will be dedicated to leveraging tools from integral quadratic constraints (IQCs) to construct accelerated EM-like algorithms by including artificial control inputs optimally designed through feedback.

\addtolength{\textheight}{-12cm}   

\bibliographystyle{IEEEtran}
\bibliography{IEEEabrv,mybibfile}

\end{document}